# Two Results for the Omega Limit Sets of Dynamical Systems

**Iasson Karafyllis**

*Dept. of Mathematics, National Technical University of Athens,
Zografou Campus, 15780, Athens, Greece.
emails: iasonkar@central.ntua.gr , iasonkaraf@gmail.com

**Abstract**
This paper provides two results for the omega limit sets of a dynamical system. We show that omega limit sets can be estimated by using functions that satisfy different -and in many cases less demanding- assumptions than the usual assumptions in Lyapunov theorems and LaSalle's theorem.

## 1. Introduction

The main results that can allow the estimation of the $\omega-$limit sets for a dynamical system are the Lyapunov theorems and LaSalle's theorem (see [1, 5]). Various extensions of the Lyapunov theorems have been proposed in the literature (e.g., [1, 3, 6]). An extension of LaSalle's theorem was proposed in [2]. However, the problem of the estimation of the $\omega-$limit sets for a nonlinear dynamical system remains a difficult problem which does not allow accurate estimations based on arbitrary functions.

The present paper provides two results for the $\omega-$limit sets of a dynamical system. We show that $\omega-$limit sets can be estimated by using functions that satisfy different -and in many cases less demanding- assumptions than the Lyapunov theorems and LaSalle's theorem.

To the best of the author's knowledge the results are new. Since the results are rather fundamental, the author of the present paper believes that they must have been discovered by other researchers in the past. However, no mention of such results is made to all papers/books that the author knows. If the reader of the paper knows similar results in the literature then the author of the present paper would like to ask the reader to communicate with him in order to withdraw the paper.

## 2. Notation and Notions for Dynamical Systems

We consider dynamical systems of the form

$$\dot{x} = f(x) \quad , \quad x \in \mathbb{R}^n \tag{1}$$

where $f : \mathbb{R}^n \to \mathbb{R}^n$ is a locally Lipschitz vector field. We use the notation $\phi(t, x_0)$ for the unique solution $x(t)$ at time $t \geq 0$ of the initial-value problem (1) with initial condition $x(0) = x_0 \in \mathbb{R}^n$.

Standard results in differential equations guarantee that for every $x_0 \in \mathbb{R}^n$ there exists $t_{\max}(x_0) \in (0, +\infty]$ such that the solution $\phi(t, x_0)$ is defined for $t \in [0, t_{\max}(x_0))$. Furthermore, if $t_{\max}(x_0) < +\infty$ then $\limsup_{t \to t_{\max}^-}(|\phi(t, x_0)|) = +\infty$.

We say that the non-empty set $A \subseteq \mathbb{R}^n$ is positively invariant for (1) if for every $x_0 \in A$ it holds that $\phi(t, x_0) \in A$ for all $t \in [0, t_{\max}(x_0))$. We say that the non-empty set $A \subseteq \mathbb{R}^n$ is negatively invariant for (1) if the following implication holds: if $\phi(t, x_0) \in A$ for certain $x_0 \in \mathbb{R}^n$ and $t \in [0, t_{\max}(x_0))$ then $x_0 \in A$. We say that the set $A \subseteq \mathbb{R}^n$ is invariant for (1) if $A \subseteq \mathbb{R}^n$ is positively and negatively invariant.

The union and the intersection of any number of (positively, negatively) invariants sets for (1) is also a (positively, negatively) invariant set for (1) (provided that the intersection is not empty). For a set $A \subseteq \mathbb{R}^n$ that contains at least one invariant set for (1), we can define the largest invariant set for (1) in $A \subseteq \mathbb{R}^n$ to be the union of all invariant sets for (1) contained in $A \subseteq \mathbb{R}^n$.

When $A \subseteq \mathbb{R}^n$ is negatively invariant for (1) and $\Omega \subseteq \mathbb{R}^n$ is positively invariant for (1) with $\Omega \setminus A \neq \emptyset$ then $\Omega \setminus A$ is positively invariant.

For every $x_0 \in \mathbb{R}^n$ we define $\omega(x_0) \subseteq \mathbb{R}^n$ to be the $\omega-$limit set of $x_0 \in \mathbb{R}^n$ using the formula:

$$\omega(x_0) = \left\{ y \in \mathbb{R}^n : \exists \{t_k \geq 0 : k = 1, 2, \ldots\} \text{ with } \lim_{k \to +\infty}(t_k) = +\infty \text{ and } \lim_{k \to +\infty}(\phi(t_k, x_0)) = y \right\} \quad (2)$$

We also define the positive orbit of $x_0 \in \mathbb{R}^n$ by the formula:

$$\gamma^+(x_0) = \left\{ \phi(t, x_0) : t \in [0, t_{\max}(x_0)) \right\} \quad (3)$$

Clearly, there is no guarantee that $\omega(x_0) \subseteq \mathbb{R}^n$ is non-empty (even if the solution $\phi(t, x_0)$ is defined for $t_{\max}(x_0) = +\infty$). When $\gamma^+(x_0) \subseteq \mathbb{R}^n$ is bounded then $\omega(x_0) \subseteq \mathbb{R}^n$ is a non-empty, compact, invariant and connected set (see [5]).

We use the notation $dist(x, A)$ for the distance of a point $x \in \mathbb{R}^n$ from a non-empty set $A \subseteq \mathbb{R}^n$, i.e., $dist(x, A) = \inf\{|y - x| : y \in A\}$. When $\gamma^+(x_0) \subseteq \mathbb{R}^n$ is bounded for some $x_0 \in \mathbb{R}^n$ then $\omega(x_0) \subseteq \mathbb{R}^n$ satisfies $\lim_{t \to +\infty}(dist(\phi(t, x_0), \omega(x_0))) = 0$.

For every non-empty set $A \subseteq \mathbb{R}^n$ we denote its closure by $\bar{A} \subseteq \mathbb{R}^n$, i.e., $\bar{A} = \{x \in \mathbb{R}^n : dist(x, A) = 0\}$. When $\Omega \subseteq \mathbb{R}^n$ is positively invariant for (1) then $\omega(x_0) \subseteq \bar{\Omega}$ for all $x_0 \in \Omega$.

Let $D \subseteq \mathbb{R}^n$ be a non-empty open set and let $S \subseteq \mathbb{R}^n$ be a set that satisfies $D \subseteq S \subseteq \bar{D}$. By $C^0(S;\Omega)$, we denote the class of continuous functions on $S$, which take values in $\Omega \subseteq \mathbb{R}^m$. By $C^k(S;\Omega)$, where $k \geq 1$ is an integer, we denote the class of functions on $S \subseteq \mathbb{R}^n$, which take values in $\Omega \subseteq \mathbb{R}^m$ and have continuous derivatives of order $k$. In other words, the functions of class $C^k(S;\Omega)$ are the functions which have continuous derivatives of order $k$ in $D$ that can be continued continuously to all points in $S$. When $\Omega = \mathbb{R}$ then we write $C^0(S)$ or $C^k(S)$.

## 3. Main Results

Our first main result provides an estimation for the $\omega-$limit sets of a dynamical system and can be applied to cases where LaSalle's theorem cannot provide accurate estimates of the $\omega-$limit sets.

**Theorem 1:** *Consider system (1) and let $\Omega \subseteq \mathbb{R}^n$ be a positively invariant compact set for (1). Suppose that there exists a non-empty set $A \subseteq \mathbb{R}^n$ and a function $V \in C^1(S)$ with $\Omega \setminus A \subseteq S$ such that the following properties hold:*

**(i)** $\nabla V(x) f(x) \neq 0$ *for all* $x \in \Omega \setminus A$,

**(ii)** *the set $A \subseteq \mathbb{R}^n$ is negatively invariant for (1).*

*Let $M \subseteq \bar{A} \cap \Omega$ be the largest invariant set in $\bar{A} \cap \Omega$. Then the inclusion $\omega(x_0) \subseteq M$ holds for every $x_0 \in \Omega$.*

The following example shows that Theorem 1 can be applied in cases where LaSalle's theorem cannot give accurate estimates of the $\omega-$limit sets.

**Example:** This example was first studied in [4] where it was shown that non-uniform attractivity properties have to be used. Consider the planar system

$$\begin{aligned} \dot{x}_1 &= -|x_2|x_2 \\ \dot{x}_2 &= |x_2|x_1 \\ x &= (x_1, x_2)' \in \mathbb{R}^2 \end{aligned} \quad (4)$$

The $\omega-$limit set can be found for every initial condition:

$$\text{If } x_2(0) = 0 \text{ then } \omega(x(0)) = \{(x_1(0), 0)\}$$

$$\text{If } x_2(0) > 0 \text{ then } \omega(x(0)) = \left\{\left(-\sqrt{x_1^2(0) + x_2^2(0)}, 0\right)\right\}$$

If $x_2(0) < 0$ then $\omega(x(0)) = \left\{\left(\sqrt{x_1^2(0) + x_2^2(0)}, 0\right)\right\}$

This can be seen by using polar coordinates $x_1 = r\cos(\theta)$, $x_2 = r\sin(\theta)$ for which we get $\dot{r} = 0$ and $\dot{\theta} = r|\sin(\theta)| \geq 0$. If we were to apply LaSalle's theorem with a function $V \in C^1(S)$ with $S \subseteq \mathbb{R}^2$ being a set that contains a ball $B$ centered at $0 \in \mathbb{R}^2$, then LaSalle's theorem would require that $V(\phi(t, x_0))$ is non-increasing along the solutions of (4) that start in $B \subseteq \mathbb{R}^2$. Let $\rho > 0$ be less than the radius of $B \subseteq \mathbb{R}^2$. Taking into account the solutions of (4), we get:

$$V(-\rho, 0) \leq V(\rho\cos(\theta), \rho\sin(\theta)), \text{ for all } \theta \in (0, \pi] \tag{5}$$

$$V(\rho, 0) \leq V(\rho\cos(\theta), \rho\sin(\theta)), \text{ for all } \theta \in (\pi, 2\pi] \tag{6}$$

Exploiting (5), (6) and continuity of $V$, we get:

$$V(-\rho, 0) \leq V(\rho\cos(\theta), \rho\sin(\theta)) \leq V(\rho, 0), \text{ for all } \theta \in [0, \pi] \tag{7}$$

$$V(\rho, 0) \leq V(\rho\cos(\theta), \rho\sin(\theta)) \leq V(-\rho, 0), \text{ for all } \theta \in [\pi, 2\pi] \tag{8}$$

It follows from (7) and (8) that $V$ has to be a function of only one variable, namely $V(x_1, x_2) = g(x_1^2 + x_2^2)$. Consequently, $\nabla V(x) f(x) = 0$ for (4) for all $x \in \mathbb{R}^2$ in $B \subseteq \mathbb{R}^2$. Since any circle $x_1^2 + x_2^2 = \rho^2$ is invariant for (4), the only piece of information that we can gain from LaSalle's theorem is that $\omega(x(0)) \subseteq \{x \in \mathbb{R}^2 : x_1^2 + x_2^2 = x_1^2(0) + x_2^2(0)\}$.

On the other hand, we notice that the set $A = \{(x_1, 0) \in \mathbb{R}^2 : x_1 \in \mathbb{R}\}$ is invariant for (4). Moreover, the function $V(x) = x_1$ is a function that satisfies all assumptions of Theorem 1 (notice that $\dot{V}(x) = -|x_2|x_2 \neq 0$ for all $x \in \mathbb{R}^2 \setminus A$) with $\Omega = \{(x_1, x_2) \in \mathbb{R}^2 : x_1^2 + x_2^2 = r^2\}$ and arbitrary $r \geq 0$. Therefore, Theorem 1 gives $\omega(x(0)) \subseteq \{(-r, 0), (r, 0)\}$ for every $x(0) \in \Omega$. ◁

Our second main result provides estimates for the $\omega$-limit sets of (1) based on an *arbitrary* function $V \in C^1(\mathbb{R}^n)$.

**Theorem 2:** *Consider system (1) and suppose that the positive orbit $\gamma^+(x_0) = \{\phi(t, x_0) : t \geq 0\}$ from the point $x_0 \in \mathbb{R}^n$ is bounded. Let the function $V \in C^1(\mathbb{R}^n)$ be given and define the set $S = \{x \in \mathbb{R}^n : \nabla V(x) f(x) = 0\}$. Then $\omega(x_0) \cap S \neq \emptyset$.*

*Furthermore, define the set $A = \{x \in \mathbb{R}^n : \exists t \geq 0 \text{ with } \phi(t, x) \in S\}$. Let $M \subseteq \overline{A}$ be the largest invariant set in $\overline{A}$. Then the inclusion $\omega(x_0) \subseteq M$ holds.*

**Remark:** It should be noticed that the assumptions of Theorem 2 guarantee that the set $S = \{x \in \mathbb{R}^n : \nabla V(x) f(x) = 0\}$ is non-empty for every $V \in C^1(\mathbb{R}^n)$ (see the proof of Theorem 2 below). It should also be noticed that if we use the terminology in [1] and if for every $x_0 \in \mathbb{R}^n$ the positive orbit $\gamma^+(x_0) = \{\phi(t, x_0) : t \geq 0\}$ is bounded then Theorem 2 guarantees that for every $V \in C^1(\mathbb{R}^n)$ the set $S = \{x \in \mathbb{R}^n : \nabla V(x) f(x) = 0\}$ is a global weak attractor. If the set $S = \{x \in \mathbb{R}^n : \nabla V(x) f(x) = 0\}$ is negatively invariant for (1) for certain $V \in C^1(\mathbb{R}^n)$ and if for every $x_0 \in \mathbb{R}^n$ the positive orbit $\gamma^+(x_0) = \{\phi(t, x_0) : t \geq 0\}$ is bounded then Theorem 2 guarantees that the set $S = \{x \in \mathbb{R}^n : \nabla V(x) f(x) = 0\}$ is a global attractor.

## 4. Proofs

**Proof of Theorem 1:** The conclusion of the theorem is trivial when $\Omega \setminus A = \emptyset$. Therefore, we assume next that $\Omega \setminus A \neq \emptyset$.

It suffices to prove that the inclusion $\omega(x_0) \subseteq M$ holds for every $x_0 \in \Omega \setminus A$. Indeed, if the inclusion $\omega(x_0) \subseteq M$ holds for every $x_0 \in \Omega \setminus A$ then the conclusion of the theorem follows, i.e., the inclusion $\omega(x_0) \subseteq M$ holds for every $x_0 \in \Omega$. To see this, let arbitrary $x_0 \in \Omega \cap A$ be given. By compactness and positive invariance of $\Omega \subseteq \mathbb{R}^n$, the limit set $\omega(x_0)$ is non-empty with $\omega(x_0) \subseteq \Omega$. We have the following cases.

<u>Case I:</u> There exists $T > 0$ such that $\phi(T, x_0) \notin A$. Thus, positive invariance of $\Omega \subseteq \mathbb{R}^n$ implies that $\phi(T, x_0) \in \Omega \setminus A$. Since $\omega(\phi(T, x_0)) = \omega(x_0)$ and $\omega(\phi(T, x_0)) \subseteq M$, it follows that $\omega(x_0) \subseteq M$ holds.

<u>Case II:</u> $\phi(t, x_0) \in A$ for all $t \geq 0$. Consequently, in this case $\omega(x_0) \subseteq \overline{A}$. Since $\omega(x_0)$ is invariant for (1), it follows that the inclusion $\omega(x_0) \subseteq M$ holds.

Therefore, we next prove that the inclusion $\omega(x_0) \subseteq M$ holds for every $x_0 \in \Omega \setminus A$.

Let arbitrary $x_0 \in \Omega \setminus A$ be given. By compactness and positive invariance of $\Omega \subseteq \mathbb{R}^n$, the limit set $\omega(x_0)$ is non-empty with $\omega(x_0) \subseteq \Omega$.

Notice that by virtue of assumption (ii) the set $\Omega \setminus A$ is positively invariant for (1) and $\phi(t, x_0) \in \Omega \setminus A$ holds for all $t \geq 0$. Define the function

$$p(t) = \dot{V}(\phi(t, x_0)) = \nabla V(\phi(t, x_0)) f(\phi(t, x_0)) \tag{9}$$

Notice that assumption (i) implies that $p(t) \neq 0$ for all $t \geq 0$. Without loss of generality, we may assume that $p(0) = \nabla V(x_0) f(x_0) > 0$ (everything that follows can also be done if $p(0) = \nabla V(x_0) f(x_0) < 0$).

We claim that $p(t) > 0$ for all $t \geq 0$. Indeed, suppose that there exists $T > 0$ such that $p(T) \leq 0$. Then, since $p(t) \neq 0$ for all $t \geq 0$, we get that $p(T) < 0$. By continuity of the function $\mathbb{R}_+ \ni t \to p(t)$ and the facts that $p(0) > 0$, $p(T) < 0$ there exists $\tau \in (0, T)$ such that $p(\tau) = 0$; a contradiction with the fact that $p(t) \neq 0$ for all $t \geq 0$.

Consequently, definition (9) guarantees that the mapping $\mathbb{R}_+ \ni t \mapsto V(\phi(t, x_0))$ is non-decreasing. Therefore, either $\lim_{t \to +\infty} (V(\phi(t, x_0))) = +\infty$ or $\lim_{t \to +\infty} (V(\phi(t, x_0))) \in \mathbb{R}$.

We distinguish the following cases:

<u>Case 1:</u> $\lim_{t \to +\infty} (V(\phi(t, x_0))) = +\infty$.

We claim that $\omega(x_0) \subseteq A$.

Indeed, suppose that there exists $y \in \omega(x_0) \subseteq \Omega$ with $y \notin A$. Since $y \in \omega(x_0)$ there exists a sequence $\{t_i \geq 0 : i = 0, 1, ...\}$ with $\lim_{i \to +\infty} (\phi(t_i, x_0)) = y$ and $\lim_{i \to +\infty} (t_i) = +\infty$. Consequently, $\lim_{i \to +\infty} (V(\phi(t_i, x_0))) = +\infty$. Continuity of $V$ on $\Omega \setminus A$, the fact that $\phi(t_i, x_0) \in \Omega \setminus A$ holds for all $i \geq 0$ and the fact that $y \in \Omega \setminus A$ imply that $\lim_{i \to +\infty} (V(\phi(t_i, x_0))) = V(y) < +\infty$; a contradiction.

<u>Case 2:</u> $\lim_{t \to +\infty} (V(\phi(t, x_0))) = L \in \mathbb{R}$.

We claim that $\omega(x_0) \subseteq A$.

Indeed, suppose that $\omega(x_0) \setminus A \neq \emptyset$. Let arbitrary $y \in \omega(x_0) \setminus A$ be given. Then there exists a sequence $\{t_i \geq 0 : i = 0, 1, ...\}$ with $\lim_{i \to +\infty} (\phi(t_i, x_0)) = y$ and $\lim_{i \to +\infty} (t_i) = +\infty$. Consequently, $\lim_{i \to +\infty} (V(\phi(t_i, x_0))) = L$. Continuity of $V$ on $\Omega \setminus A$, the fact that $\phi(t_i, x_0) \in \Omega \setminus A$ holds for all $i \geq 0$ and the fact that $y \in \Omega \setminus A$ imply that $\lim_{i \to +\infty} (V(\phi(t_i, x_0))) = V(y)$. Therefore, $V(y) = L$.

Since $y \in \omega(x_0) \setminus A$ is arbitrary, we conclude that $V(y) = L$ for all $y \in \omega(x_0) \setminus A$.

Let $z \in \omega(x_0) \setminus A$ be given. Since $\omega(x_0)$ is invariant for (1), it follows that $\phi(t, z) \in \omega(x_0)$ for all $t \geq 0$. Since the set $\Omega \setminus A$ is positively invariant for (1), it follows that $\phi(t, z) \in \Omega \setminus A$ for all $t \geq 0$. Thus, $\phi(t, z) \in \omega(x_0) \setminus A$ for all $t \geq 0$. Consequently, $V(\phi(t, z)) = L$ for all $t \geq 0$. This implies $\dot{V}(\phi(t,z)) = 0$ for all $t \geq 0$, which gives $\dot{V}(z) = \nabla V(z) f(z) = 0$; a contradiction with assumption (i).

Therefore, we have shown that the inclusion $\omega(x_0) \subseteq A$ holds. Since $\omega(x_0) \subseteq \Omega$ is invariant for (1), it follows that the inclusion $\omega(x_0) \subseteq M$ holds.

Since $x_0 \in \Omega \setminus A$ is arbitrary, it follows that the inclusion $\omega(x_0) \subseteq M$ holds if $x_0 \in \Omega \setminus A$. The proof is complete. ◁

**Proof of Theorem 2:** First of all, we notice that the assumptions of the theorem guarantee that the set $S = \{x \in \mathbb{R}^n : \nabla V(x) f(x) = 0\}$ is non-empty. Indeed, if the set $S = \{x \in \mathbb{R}^n : \nabla V(x) f(x) = 0\}$ were empty then we would have either $\nabla V(x) f(x) > 0$ for all $x \in \mathbb{R}^n$ or $\nabla V(x) f(x) < 0$ for all $x \in \mathbb{R}^n$. If $\nabla V(x) f(x) > 0$ for all $x \in \mathbb{R}^n$ and if $\Omega \subseteq \mathbb{R}^n$ is a compact set that contains the bounded set $\gamma^+(x_0) = \{\phi(t, x_0) : t \geq 0\}$ then we would be able to define $b = \min\{\nabla V(x) f(x) : x \in \Omega\} > 0$ and $C = \max\{V(x) : x \in \Omega\}$. We would then have for all $t \geq 0$:

$$C \geq V(\phi(t, x_0)) = V(x_0) + \int_0^t \nabla V(\phi(\tau, x_0)) f(\phi(\tau, x_0)) d\tau \geq V(x_0) + bt \tag{10}$$

However, since $b > 0$, inequality (10) cannot hold for all $t \geq 0$; a contradiction. A similar contradiction is obtained in the case that $\nabla V(x) f(x) < 0$ for all $x \in \mathbb{R}^n$.

The proof of the theorem is made by means of a contradiction argument.

Suppose that $\omega(x_0) \cap S = \emptyset$.

Then for every sequence $\{t_k \geq 0 : k = 1, 2, \ldots\}$ with $\lim_{k \to +\infty}(t_k) = +\infty$ we cannot have $\lim_{k \to +\infty}(dist(\phi(t_k, x_0), S)) = 0$. We prove this fact by contradiction. Suppose that there exists a sequence $\{t_k \geq 0 : k = 1, 2, \ldots\}$ with $\lim_{k \to +\infty}(t_k) = +\infty$ and $\lim_{k \to +\infty}(dist(\phi(t_k, x_0), S)) = 0$. Notice that since the set $\{\phi(t, x_0) : t \geq 0\}$ is bounded, there exists a subsequence of $\{t_k \geq 0 : k = 1, 2, \ldots\}$, still denoted by $\{t_k \geq 0 : k = 1, 2, \ldots\}$ with $\lim_{k \to +\infty}(t_k) = +\infty$, $\lim_{k \to +\infty}(dist(\phi(t_k, x_0), S)) = 0$ and $\lim_{k \to +\infty}(\phi(t_k, x_0)) = x^* \in \omega(x_0)$. Clearly, continuity of the mapping $x \mapsto dist(x, S)$ implies that $dist(x^*, S) = 0$. Since $S = \{x \in \mathbb{R}^n : \nabla V(x) f(x) = 0\}$ is closed, this implies that $x^* \in S$ which shows that $x^* \in \omega(x_0) \cap S = \emptyset$; a contradiction.

We next claim that $\liminf_{t \to +\infty}(dist(\phi(t, x_0), S)) > 0$. Indeed, notice that the fact that the set $\{\phi(t, x_0) : t \geq 0\}$ is bounded and the fact that $dist(\phi(t, x_0), S) \geq 0$, imply that

$$0 \le \liminf_{t \to +\infty} \left( dist(\phi(t, x_0), S) \right) \le \limsup_{t \to +\infty} \left( dist(\phi(t, x_0), S) \right) < +\infty$$

Consequently, in order to prove the claim $\liminf_{t \to +\infty} \left( dist(\phi(t, x_0), S) \right) > 0$ it suffices to show that $\liminf_{t \to +\infty} \left( dist(\phi(t, x_0), S) \right) \ne 0$. The proof of the claim is made by contradiction. Suppose that $\liminf_{t \to +\infty} \left( dist(\phi(t, x_0), S) \right) = 0$. Then we would have a sequence $\{ t_k \ge 0 : k = 1, 2, \ldots \}$ with $\lim_{k \to +\infty} (t_k) = +\infty$ and $\lim_{k \to +\infty} \left( dist(\phi(t_k, x_0), S) \right) = 0$; a contradiction.

Since $L = \liminf_{t \to +\infty} \left( dist(\phi(t, x_0), S) \right) > 0$, there exists $T > 0$ such that $dist(\phi(t, x_0), S) > \frac{L}{2}$ for all $t > T$. Consequently, definition $S = \{ x \in \mathbb{R}^n : \nabla V(x) f(x) = 0 \}$ implies that $\nabla V(\phi(t, x_0)) f(\phi(t, x_0)) \ne 0$ for all $t > T$. Without loss of generality, we may assume that $\nabla V(\phi(T+1, x_0)) f(\phi(T+1, x_0)) > 0$ (everything that follows can also be done if $\nabla V(\phi(T+1, x_0)) f(\phi(T+1, x_0)) < 0$).

We claim that $\nabla V(\phi(t, x_0)) f(\phi(t, x_0)) > 0$ for all $t \ge T + 1$. Indeed, if $\nabla V(\phi(t, x_0)) f(\phi(t, x_0)) \le 0$ for certain $t > T + 1$, then by continuity of the mapping $\tau \mapsto \nabla V(\phi(\tau, x_0)) f(\phi(\tau, x_0))$, we would obtain the existence of $\tau \in (T+1, t]$ with $\nabla V(\phi(\tau, x_0)) f(\phi(\tau, x_0)) = 0$. This contradicts the fact that $\nabla V(\phi(\tau, x_0)) f(\phi(\tau, x_0)) \ne 0$ for all $\tau > T$.

Consequently, the mapping $[T+1, +\infty) \ni t \mapsto V(\phi(t, x_0))$ is non-decreasing. Since the set $\{ \phi(t, x_0) : t \ge 0 \}$ is bounded, we also obtain that the mapping $[T+1, +\infty) \ni t \mapsto V(\phi(t, x_0))$ is bounded. Therefore $\lim_{t \to +\infty} \left( V(\phi(t, x_0)) \right) = G \in \mathbb{R}$.

Let arbitrary $y \in \omega(x_0)$ be given. Then there exists a sequence $\{ t_i \ge 0 : i = 0, 1, \ldots \}$ with $\lim_{i \to +\infty} (\phi(t_i, x_0)) = y$ and $\lim_{i \to +\infty} (t_i) = +\infty$. Consequently, $\lim_{i \to +\infty} \left( V(\phi(t_i, x_0)) \right) = G$. Continuity of $V$ implies that $\lim_{i \to +\infty} \left( V(\phi(t_i, x_0)) \right) = V(y)$. Therefore, $V(y) = G$.

Since $y \in \omega(x_0)$ is arbitrary, we conclude that $V(y) = G$ for all $y \in \omega(x_0)$.

Let $z \in \omega(x_0)$ be given. Since $\omega(x_0)$ is invariant for (1), it follows that $\phi(t, z) \in \omega(x_0)$ for all $t \ge 0$. Consequently, $V(\phi(t, z)) = G$ for all $t \ge 0$. This implies $\dot{V}(\phi(t, z)) = 0$ for all $t \ge 0$, which gives $\dot{V}(z) = \nabla V(z) f(z) = 0$. Therefore, definition $S = \{ x \in \mathbb{R}^n : \nabla V(x) f(x) = 0 \}$ implies that $z \in S$. This implies that $z \in \omega(x_0) \cap S = \emptyset$; a contradiction.

Notice that the set $A = \{x \in \mathbb{R}^n : \exists t \geq 0, \phi(t,x) \in S\}$ is negatively invariant for (1) and satisfies $S \subseteq A$. Consequently, definition $S = \{x \in \mathbb{R}^n : \nabla V(x) f(x) = 0\}$ implies that $\nabla V(x) f(x) \neq 0$ for all $x \in \mathbb{R}^n \setminus A$. Applying Theorem 1 with $\Omega \subseteq \mathbb{R}^+$ being the closure of the positive orbit $\{\phi(t, x_0) : t \geq 0\}$, we conclude that the inclusion $\omega(x_0) \subseteq M$ holds. The proof is complete. ◁